\documentstyle[12pt]{article}
\begin{document}

\title{The Riemann problem with additional
singularities}
\author{M. Kudryavtsev}

\maketitle
\begin{center}
{\it The Institute for Low Temperature Physics and Engineering\\
of the National Academy of Science of Ukraine,\\
\smallskip
Ukraine, Kharkov, 61164, Lenine Ave., 47}\\
\medskip
E-mail: kudryavtsev@ilt.kharkov.ua\\
\medskip
\end{center}

\abstract{

   The Riemann problem is studied in the case when the unknown function
has nonisolated singularities, concentrated on the real axis.
The problem is used for the factorization of functions,
holomorphic outside of the unit circle and the real axis, in the form
of product of two functions which have singularities on the given set
of the real axis.
}

\bigskip \bigskip \bigskip

The classical Riemann problem with zeroes is well known (see, e.g. [1]). It
is required to construct two functions $\psi _{1}(z)$, $\psi _{2}(z)$, such
that $\psi _{1}(z)$ is analytic inside the closed contour $\Gamma $ and has
inside the contour $n$ zeroes $\lambda _{1},\ldots ,\lambda _{n}$, and ${%
\frac{1}{\psi _{2}(z)}}$ is analytic outside $\Gamma $ and has outside the
contour the zeroes $\mu _{1},\ldots ,\mu _{n}$. In addition, on the contour $%
\Gamma $ the following relation is required:
$$
\psi _{1}(\xi )=G(\xi )\psi _{2}(\xi ),\quad \xi \in \Gamma ,\eqno(1)
$$
where $G(\xi )$ is a given complex-valued function on the contour.

In this work the Riemann problem with additional singularities is proposed
and solved in the case when the contour $\Gamma $ is the unit circle and the
zeroes and the singularities of the functions $\psi _{1}(z)$ and
$\psi_{2}(z)$, including non isolated singularities, are concentrated
on the real axis. A particular case of such a problem is used in [2].

\bigskip

{\it Notation 1.} \ If the function $f(z)$ is holomorphic everywhere in
complex plane outside the unite circle ${\bf T}=\{\xi ,\,|\xi |=1\}$ and the
real axis ${\bf R}$, then the superindices $^{+}$ and $^{-}$ denote the
limit value of the function from inside ($^{+}$) and from the outside ($^{-}$%
) of the unit circle or from above ($^{+}$) and from below ($^{-}$) of the
real axis (assuming that these limits exist):
$$
   f^{\pm }(\xi )=\lim_{\varepsilon \rightarrow +0}f((1\mp \varepsilon )\xi),
    \quad |\xi |=1,
$$
$$
   f^{\pm }(x)=\lim_{\varepsilon \rightarrow +0}f(x\pm i\varepsilon ),
   \quad -\infty <x<\infty .
$$

\medskip

Now we formulate the Riemann problem with additional singularities. Let the
function
$$
   G(e^{i\theta })=\mu (\theta )e^{i\hat{\varphi}(\theta )},
   \quad -\pi <\theta \leq \pi ,
$$
be given on the unit circle, where $\ln \mu (\theta )$ and
$\hat{\varphi}(\theta )$ are summable functions ($\mu (\theta )>0$,
$\hat{\varphi}(\theta)=\overline{\hat{\varphi}(\theta )}$), and
$\hat{\varphi}(\theta )$ is bounded. Let $\varphi (t)=\overline{\varphi (t)}$
be a bounded summable function on the real axis, which vanishes at least on
one interval $\Delta^{\prime }\subset (-1,1)$ and on one interval
$\Delta^{\prime \prime}\subset (-\infty ,-1)\cup (1,\infty )$.
It is required to construct a holomorphic function $R(z)$, which do not have
zeroes outside the unit circle and the real axis, and such that
$$
  {\frac{R^{+}(e^{i\theta })}{R^{-}(e^{i\theta })}}
   =\mu (\theta )e^{i\hat{\varphi}(\theta )},
   \quad -\pi <\theta \leq \pi ,\eqno(2)
$$
$$
   {\frac{R^{+}(t)}{R^{-}(t)}}=e^{i\varphi (t)},\quad -\infty <t<\infty ,
    \eqno (3)
$$
assuming that these limits exist. Here we require the exact equality of the
arguments in (2) and (3): $\arg {\frac{R^{+}(e^{i\theta })}{R^{-}(e^{i\theta
})}}=\hat{\varphi}(\theta )$,
$\arg {\frac{R^{+}(t)}{R^{-}(t)}}=\varphi (t)$, where
$\arg R^{\pm }(\xi ),\ |\xi |=1,$ and
$\arg R^{\pm }(t),-\infty<t<\infty $, are defined in the following way.
Let us fix points $t^{\prime}\in \Delta ^{\prime }$,
$t^{\prime \prime }\in \Delta ^{\prime \prime }$,
in which $R(z)$ is holomorphic and let $-\pi <\arg R(t^{\prime })\leq \pi $,
$-\pi <\arg R(t^{\prime \prime })\leq \pi $. According to the assumption,
the function $R(z)$ has four connected components of holomorphy in which it
does not vanish. It is the parts of the upper and lower halfplane restricted
by the unit circle. Connecting the point $z$
(${\rm Im}\,z\neq 0,\ |z|\neq 1$) with the point $t^{\prime }$ or
$t^{\prime \prime }$ by a continuous curve, lying in one of the four
components, we observe the continuous change of the argument of the
function $R(z)$ along this curve. Now the argument $\arg R(z)$ is defined
uniquely. As $z$ tends to the real axis (resp. to the unit circle), we find
$\arg R^{\pm }(t),-\infty <t<\infty $ (resp. $\arg R^{\pm }(\xi ),\ |\xi |=1$).

It is easy to see that the Riemann problem with additional singularities is
a generalization of the classical Riemann problem with zeroes when the
zeroes are concentrated on the real axis. In fact, if on some interval
$\Delta \in {\bf R},\ \pm 1\not\in \Delta ,$
$\varphi (t)=k\pi ,\ k\in {\bf Z}$, then this means that $R(z)$ is
holomorphic on $\Delta $. If in the right and left half-neighborhoods
of the point $t_{0}$ the function $\varphi (t)$ is constant and divisible
by $\pi $, and at the point the function has a jump of the form $k\pi $,
this means that $R(z)$ has at the point $t_{0}$ a pole ($k>0$) or a zero
($k<0$) of the order $|k|$. More complicated behavior of $\varphi (t)$
implies more complicated character of the singularities of $R(z)$.

\medskip

{\it Notation 2.} Introduce
$$ P(z,\gamma )\ =\ \exp
    \Bigl\{{\frac{1}{{\pi }}}\int\limits_{-\infty}^{\infty }\gamma (t)
      ({\frac{1}{{t-z}}}-{\frac{t}{{1+t^{2}}}})dt\Bigr\}, \eqno (4)
$$
$$
   \hat{P}(z,\hat{\gamma})\ =\
    \exp \biggl\{-{\frac{1}{{2\pi i}}}{\int\limits_{-\pi }^{\pi }}
    {\frac{{e^{i\theta }+z}}{{e^{i\theta }-z}}}\ \hat{\gamma}(\theta )
    d\theta \biggr\},
$$
where $\gamma (t)=\overline{\gamma (t)},\ -\infty <t<\infty ,$ and
$\hat{\gamma}(\theta )=\overline{\hat{\gamma}(\theta )},\ -\pi <t<\pi ,$
are bounded measurable functions. The function $P(z,\gamma )$ is defined
and holomorphic at least for non real $z$, and $\hat{P}(z,\hat{\gamma})$,
resp., is defined and holomorphic inside and outside the unit disk. The
formulae of Plemelj-Sokhotsky imply the following equalities, connecting
the limit values of $P(z,\gamma )$ and $\hat{P}(z,\hat{\gamma})$ on the
real axis and the unit circle resp.:
$$
   \arg P^{+}(t,\gamma )=-\arg P^{-}(t,\gamma )=\gamma (t),\quad
   |P^{+}(t,\gamma )|=|P^{-}(t,\gamma )|,\qquad -\infty <t<\infty ,\eqno(5)
$$
$$
    \arg \hat{P}^{+}(e^{i\theta },\hat{\gamma})
       =-\arg \hat{P}^{-}(e^{i\theta },
    \hat{\gamma})=\hat{\gamma}(\theta ),
    \quad |\hat{P}^{+}(e^{i\theta },\hat{\gamma})|
     =|\hat{P}^{-}(e^{i\theta },\hat{\gamma})|,\qquad -\pi <\theta <\pi.
  \eqno(6)
$$

\smallskip

Let us define the functions
$$
   R^{(1)}(z)=P(z,{\frac{\varphi }{2}}),
  \qquad R^{(2)}(z)=\hat{P}(z,{\frac{\hat{\varphi}}{2}}),
$$
$$
   R_{\mu }(z)=\exp \biggl\{{\frac{1}{2\pi }}\int\limits_{-\pi }^{\pi }
    {\frac{e^{i\theta }}{e^{i\theta }-z}}\ln \mu (\theta )\ d\theta \biggr\}.
   \eqno(7)
$$
It is easy to obtain from the formulae of Plemel-Sokhotsky that
$$
   {\frac{|R_{\mu }^{+}(e^{i\theta })|}{|R_{\mu }^{-}(e^{i\theta })|}}
  =\mu(\theta),  \quad
  \arg R_{\mu }^{+}(e^{i\theta })=\arg R_{\mu }^{-}(e^{i\theta}),
  \qquad -\pi <\theta <\pi .\eqno(8)
$$

Thus, the following theorem is a simple consequence of the equalities (5),
(6) and (8):

\bigskip

{\bf Theorem 1.} \ {\it The function $R(z)=R^{(1)}(z)R^{(2)}(z)R_{\mu }(z)$
is a solution of the problem (2), (3). }

\bigskip

We observe that the summability of the functions $\ln \mu (\theta )$,
$\hat{\varphi}(\theta )$, $\varphi (t)$ being provided, we have the
existence of the limits $R^{\pm }(\xi )$ and $R^{\pm }(t)$ almost
everywhere on the unit circle and on the real axis. and the functions
$\ln R((1\mp \varepsilon )\xi),\ |\xi |=1$, and
$\ln R(t\pm i\varepsilon ),\ t\in {\bf R}$, converge respectively to
$\ln R^{\pm }(\xi )$ and $\ln R^{\pm }(t)$ with respect to the metric of
$L^{1}$ (see, e.g., [3]). Without paying attention to the problems of
convergence, we will apply the Riemann problem to the factorization
of functions with singularities on the unit circle and the
real axis.

\medskip

{\it Notation 3.} We define on the real axis the map $V$ of the symmetry
with respect to the unit circle:
$$
V(t)=t^{-1},\quad t\in {\bf R}\backslash \{0\}.
$$
For a set $A\subset {\bf R}\backslash \{0\}$ and a function $\rho (t)$,
defined on ${\bf R}\backslash \{0\}$ we denote
$$
V(A)=\{t\ |\ t^{-1}\in A\},
$$
$$
V(\rho )(t)=\rho (t^{-1}).
$$

The main result of the work is the following

{\bf Theorem 2.} \ {\it Let the function $N(z)$ be holomorphic and not
vanishing outside the unit circle ${\bf T}$ and a certain closed set $\Sigma
\subset {\bf R}$ of the real axis and satisfies the following conditions:

1) In the domain of holomorphy of the function $N(z)$ we have
$$
   N(z^{-1})=N(z),\quad N(\overline{z})=\overline{N(z)},
  \qquad z\not\in \Sigma \cup {\bf T}.\eqno(9)
$$

2) Their exists such positive constant $C>0$, that in the upper
semi-disk
$$
|\arg N(z)|<C,\quad \quad |z|<1,\ {\rm Im}\,z>0.
$$

}

{\it If $\Sigma =\Omega _{1}\cup V(\Omega _{1})\cup \Omega _{2}$, where
the sets $\Omega _{1}$, $V(\Omega _{1})$, and  $\Omega _{2}=V(\Omega _{2})$
have the mutual positive distances, then there exists the function $R(z)$,
holomorphic outside the set $\Omega \equiv \Omega _{1}\cup \Omega _{2}$
and outside the unit circle {\bf T}, such that
$$
N(z)=R(z)R(z^{-1})\eqno(10),
$$
and
$$
   {\frac{|R^{+}(e^{i\theta })|}{|R^{-}(e^{i\theta })|}}=\mu (\theta ),
   \eqno(11)
$$
where $\mu (\theta )=\mu (-\theta )>0,\ -\pi <\theta <\pi $, is an
arbitrary even function with summable logarithm. }

We remark that this theorem reduces to the Riemann problem with additional
singularities. Essentially, we construct the function $R(z)$ so that the
limit values of the arguments of the function $R(z)R(z^{-1})$ equal the
limit values of the argument of $N(z)$. We also require that the function $%
R(z)$ only have singularities onto the unit circle {\bf T} and onto the set $%
\Omega $, and satisfy the additional condition (11) with practically
arbitrary $\mu (\theta )$. In order to prove the theorem, we will need two
simple lemmas demonstrating the properties of the functions $P(z,\gamma )$
and $\hat{P}(z,\hat{\gamma})$.

{\bf Lemma 1.} \ {\it The functions $P(z,\gamma)$ and $\hat P(z,\hat\gamma)$
in their domain of holomorphy satisfy the properties:
$$
   P(z,\gamma_1+\gamma_2)=P(z,\gamma_1) P(z,\gamma_2), \eqno(12)
$$
$$
   P(\overline{z},\gamma)  = \overline{P(z,\gamma)}, \eqno(13)
$$
$$
   P(z^{-1},\gamma)  = P(z,-V(\gamma)), \eqno(14)
$$
$$
   \hat P(z,\hat\gamma_1+\hat\gamma_2)  =
   \hat P(z,\hat\gamma_1) \hat P(z,\hat\gamma_2), \eqno(15)
$$
If $\hat\gamma(\theta)=-\hat\gamma(-\theta)$ is odd, then}
$$
\hat P(z^{-1},\hat\gamma) = \hat P(z,\hat\gamma)). \eqno(16)
$$

{\bf Lemma 2.} \ {\it 1. Let $f_{1}(z)$ be holomorphic function of bounded
argument in the halfplane $\{{\rm Im}\,z>0\}$, which do not have zeroes in
the halfplane. Then it can be represented in the form
$$
   f_{1}(z)=C_{1}P(z,\eta _{1}),\quad {\rm Im}\,z>0,\eqno(17)
$$
where $C_{1}$ is a positive constant and
$$
   \eta _{1}(t)=\arg f_{1}^{+}(t),\quad -\infty <t<\infty ,\eqno(18)
$$
is the limit value of its argument from above on the real axis. }

{\it 2. Let $f_{2}(z)$ be holomorphic function in the disk $|z|<1$ of
bounded argument. Then it can be represented in the form.
$$
   f_{2}(z)=C_{2}\hat{P}(z,\hat{\eta}_{2}),\quad |z|<1,
$$
where $C_{2}>0$, and
$$
  \hat{\eta}_{2}(\theta )=\arg f_{2}^{+}(e^{i\theta }),
  \quad -\pi <\theta <\pi,
$$
is its limit value from inside the unit circle.}

P r o o f \ of \thinspace\ lemmas. We will only prove the first part of
Lemma 2 (equalities (12)--(16) are obtained by direct calculation). The
function $f_{1}(z)$ is holomorphic and does not vanish in the connected
domain $\{{\rm Im}\,z>0\}$. Hence, we can define uniquely a logarithm $\ln
f_{1}(z)$, which is holomorphic in the upper halfplane function with bounded
imaginary part: $|{\rm Im}\,\Bigl(\ln f_{1}(z)\Bigl)|<C_{3}$. Hence, the
function $\,\Bigl(\ln f_{1}(z)+C_{3}i\Bigr)\,$ is the function of Nevanlinna
(i.e. has positive imaginary part in the upper halfplane) and can be
represented in the form (see, e.g., [4])
$$
\ln f_{1}(z)+iC_{3}=\alpha +\beta z+{\frac{1}{{\pi }}}\int\limits_{-\infty
}^{\infty }({\frac{1}{{t-z}}}-{\frac{t}{{1+t^{2}}}})d\rho (t),
$$
where the measure $d\rho (t)$ is defined by a nondecreasing function
$\rho(t)$ with
$$
  \rho (t_{2})-\rho (t_{1})
    =\lim_{\varepsilon \rightarrow 0}\int_{t_{1}}^{t^{2}}{\rm Im}\,
     \bigl(\ln f_{1}(t+i\varepsilon )+iC_{3}\bigr)dt,
$$
$\alpha \in {\bf R}$, $\beta \geq 0$. Since $\,{\rm Im}\,\ln f_{1}(z)\,$
is bounded, then $\beta =0$, and the measure $d\rho (t)$ is absolutely
continuous $d\rho (t)=(\arg f_{1}^{+}(t)+C_{3})dt$. Besides,
$$
   iC_{3}={\frac{1}{{\pi }}}\int\limits_{-\infty }^{\infty }({\frac{1}{{t-z}}}-{%
   \frac{t}{{1+t^{2}}}})C_{3}dt,
$$
hence,
$$
   \ln f_{1}(z)=\alpha +{\frac{1}{{\pi }}}\int\limits_{-\infty }^{\infty }
   ({\frac{1}{{t-z}}}-{\frac{t}{{1+t^{2}}}})\eta _{1}(t)dt,
$$
where $\eta _{1}(t)$ is defined by formula (18). This means that for the
function $f_{1}(z)$ the {\it multiplicative representation} (17) is
obtained. The second part of the lemma is proven analogously. \hfill
\rule{0.5em}{0.5em}\medskip

P r o o f \ of \thinspace\ Theorem 1. \ At first, we will present $N(z)$ in
the form of the product $P(z,\gamma )\hat{P}(z,\hat{\gamma})$. It follows
from (9) that $\arg N(z)$ is bounded in the part of the upper halfplane that
lies outside the unit disk. The function
$$
M(\lambda )=N(z)|_{z+z^{-1}=\lambda ,|z|>1}
$$
is holomorphic in the upper and lower haldplane, and its argument is
bounded. According to lemma 2,
$$
M(\lambda )=C P(\lambda ,\eta ),
$$
where $\eta (\tau )=\arg M^{+}(\tau ),\ -\infty <\tau <\infty ,$ is a
bounded function, and $C>0$ is a positive constant. (Henceforth we will
denote by $C$ positive constants, different for each case.) This
representation is true in the both halfplanes, because
$M(\overline{z})=\overline{M(z)}$, and because for $P(\lambda ,\eta )$
we have (13). Let $\chi _{\lbrack -2,2]}(t)$ and
$\chi _{(-\infty ,-2)\cup (2,\infty )}(t)$ be the indicators of the sets
$[-2,2]$ ${\bf R}\backslash \lbrack -2,2]$,
respectively. We define the functions
$$
   M_{0}(\lambda )=P(\lambda ,\chi _{\lbrack -2,2]}\eta _{0}),
$$
$$
   M_{1}(\lambda )=P(\lambda ,\chi _{(-\infty ,-2)\cup (2,\infty )}(t)\eta_{1}).
$$
According to property (12), in the domain of holomorphy of $M(\lambda )$
$$
   M(\lambda )=CM_{0}(\lambda )M_{1}(\lambda ),\eqno(19)
$$
where, according to (5), the functions $M_{0}(\lambda )$ and $M_{1}(\lambda)$
are holomorphic on ${\bf R}\backslash \lbrack -2,2]$ and $(-2,2)$,
respectively. Let us define in the plane of the parameter $z$ the functions
$$
   N_{0}(z)=N_{0}(z^{-1})\equiv M_{0}(z+z^{-1}),\quad |z|\neq 1,\qquad
   N_{1}(z)=N_{1}(z^{-1})\equiv M_{1}(z+z^{-1}),\quad z\notin \Sigma .
$$
It is evident that the functions $N_{0}(z)$ and $N_{1}(z)$ are holomorphic
and positive on the real line and the unit circle, respectively (except
possibly the points $\pm 1$). Taking into account (9), (19) and this
representation, we have at first outside, then inside the unit disk
$$
   N(z)=CN_{0}(z)N_{1}(z),\quad z\not\in \Sigma \cup {\bf T},\eqno(20)
$$
where the functions $N_{0}(z)$ and $N_{1}(z)$, according to Lemma 2, can be
represented in the multiplicative form
$$
N_{0}(z)=\hat{P}(z,\hat{\nu}_{0}),\qquad N_{1}(z)=P(z,\nu ),\eqno(21)
$$
with an odd function on $(\pi ,\pi )$
$$
   \hat{\nu}_{0}(\theta )=-\hat{\nu}_{0}(-\theta )\equiv \arg
   N_{0}^{+}(e^{i\theta })=\arg N^{+}(e^{i\theta })=%
   \cases{ -\eta(e^{i\theta}+ e^{-i\theta}), \quad 0\leq\theta<\pi, \cr
            \eta(e^{i\theta}+ e^{-i\theta}), \quad -\pi<\theta<0, \cr}
   \eqno(22)
$$
and
$$
\nu (t)=-\nu (t^{-1})\equiv \arg N_{1}^{+}(t)=\arg N^{+}(t)=
   \cases{ \eta(t+t^{-1}), \quad |t|>1, \cr
              -\eta(t+t^{-1}), \quad |t|<1. \cr}\eqno(23)
$$
We remark that Lemma 2 guarantees representation (21) for $N_{0}(z)$ inside
the disk (for $N_{1}(z)$ in the upper halfplane). However, the same
representation is also true outside the disk (resp. in the lower halfplane),
because of (16) and because $N_{0}(z)=N_{0}(z^{-1})$ (because $N_{1}(%
\overline{z})=\overline{N_{1}(z)}$). Thus, the problem of the factorization
of the function $N(z)=CN_{0}(z)N_{1}(z)$ is reduced to the problem of
factorization of two functions $N_{0}(z)$ and $N_{1}(z)$, represented in the
form (21). At first we factorize $N_{1}(z)$. Let
$$
   \Delta ={\bf R}\backslash ([-1,1]\cup \Sigma )=
   {\bf R}\backslash ([-1,1]\cup\Omega _{1}
     \cup V(\Omega _{1})\cup \Omega _{2})=\cup _{k}\Delta _{k}\,,
$$
where $\Delta _{k}=(\alpha _{k},\beta _{k})$ are mutually disjoint
intervals and $|\alpha _{k}|\geq 1,|\beta _{k}|\geq 1$.

It follows from the definition of $\Delta $ that the endpoints
$(\alpha_{k},\beta _{k})$ of the interval belong to one of the three disjoint
sets $\Omega _{1},V(\Omega _{1}),\Omega _{2}$. Moreover, the number of the
intervals, whose endpoints belong to different sets, is finite. In fact, if
the endpoints of the interval $(\alpha _{k},\beta _{k})$ belong to different
sets, then
$$
  \beta _{k}-\alpha _{k}\geq \min \{{\rm dist}(\Omega _{1},V(\Omega _{1})),
  {\rm dist}(\Omega _{1},\Omega _{2}),{\rm dist}(V(\Omega _{1}),\Omega _{2})\}.
$$
If there were infinitely many of such intervals, then (except the case when
the intervals are concentrated at the infinity) some of the intervals would
have arbitrarily small length, so that one of the distances
${\rm dist} (\Omega _{1},V(\Omega _{1}))$,
${\rm dist}(\Omega _{1},\Omega _{2})$ and
${\rm dist}(V(\Omega _{1}),\Omega _{2})$ would vanish, which contradicts
to the conditions of the theorem. (If these intervals were concentrated
in the infinity, then, according to $\Omega _{2}=V(\Omega _{2})$ and
$\Omega_{1}=V(\Omega _{1})$, we would have that the distance between
the different sets and zero vanish, which contradicts to the conditions
of the theorem, too.)

Let us choose among the intervals $\Delta _{k}$ such intervals, that one
of the endpoints belongs to $V(\Omega _{1})$, and the other belongs to
$\Omega_{1}$ or $\Omega _{2}$, or is equal to $\pm 1$. Let us renumerate
the intervals $\Delta _{k}$ so that $\Delta _{1},\ldots ,\Delta _{k_{0}}$
are the chosen intervals (their number is finite). We divide the rest of
the intervals in two groups: the intervals $\Delta _{k}^{\prime }$, whose
both endpoints belong to $V(\Omega _{1})$, and
$\Delta _{k}^{\prime \prime }$ (all the rest). Let
$$
   \alpha _{k}^{\ast }=\cases{\alpha_k, {\ \rm when\ }
                              \alpha_k\not\in V(\Omega_1),\cr
                     \alpha_k^{-1},
                              {\ \rm when\ } \alpha_k\in V(\Omega_1),\cr}
     \qquad \beta _{k}^{\ast }=
     \cases{\beta_k, {\ \rm when\ }
                                   \beta_k\not\in V(\Omega_1),\cr
                          \beta_k^{-1},
                               {\ \rm when\ } \beta_k\in V(\Omega_1).\cr}
      \qquad k=1,\ldots ,k_{0}.\eqno(24)
$$
We observe that $\{\alpha _{k}^{\ast }\}_{k=1}^{k_{0}}$,
$\{\beta _{k}^{\ast}\}_{k=1}^{k_{0}}$ and the intervals
$V(\Delta _{k}^{\prime })$, $\Delta_{k}^{\prime \prime }$ have positive
distance with the set $V(\Omega _{1})$. We denote by $\chi_1(t)$,
$\chi_2(t)$, $\chi^{(k_0)}_\Delta(t)$, $\chi^{\prime}_\Delta(t)$,
$\chi^{\prime\prime}_\Delta(t)$ the indicators of the set $\Omega_1$,
$\Omega_2$, $\cup_{k=1}^{k_0} \Delta_k$, $\cup_k \Delta^{\prime}_k$,
$\cup_k \Delta^{\prime\prime}_k$. Evidently,
$$
\begin{array}{l}
   \chi_1(t) + V(\chi_1(t)) + \chi_2(t) + \chi^{(k_0)}_\Delta(t) +
   V(\chi^{(k_0)}_\Delta)(t) + \\ +\chi^{\prime}_\Delta(t) +
   V(\chi^{\prime}_\Delta)(t) + \chi^{\prime\prime}_\Delta(t) +
   V(\chi^{\prime\prime}_\Delta)(t) = 1, \quad t \ne \pm 1.
\end{array}
\eqno(25)
$$

\medskip

{\bf Lemma 3 } \ (factorization of $N_{1}(z)$). \quad {\it The function
$N_{1}(z)$ can be factored out as follows:
$$
   N_{1}(z)\ =\ CR_{012}(z)R_{012}(z^{-1}),\eqno(26)
$$
$$
   R_{012}(z)\ =\
   R_{0}(z)R_{0}^{\prime }(z)R_{0}^{\prime \prime}(z)R_{1}(z)R_{2}(z),
$$
where
$$
   R_{0}(z)\ =\ \prod_{k=1}^{k_{0}}\Bigl({\frac{{z-\beta _{k}^{\ast }}}
   {{z-\alpha _{k}^{\ast }}}}\Bigr)^{n(k)}\,,\eqno(27)
$$
$$
   n(k)={\frac{1}{\pi }}\arg N^{+}(t)\in {\bf Z},\quad t\in \Delta _{k},\
   k=1,\ldots ,k_{0},
$$
$$
   R_{0}^{\prime }(z)\ =\ P(z,V(\chi _{\Delta }^{\prime })\nu ),  \qquad
   R_{0}^{\prime \prime }(z)=P(z,\chi _{\Delta }^{\prime \prime }\nu )\,,
   \eqno(28)
$$
$$
   R_{1}(z)\ =\ P(z,\chi _{1}\nu ),\qquad
   R_{2}(z)=P(z,{\frac{1}{2}}\chi_{2}\nu ),\eqno(29)
$$
the constant $C>0$, and the numbers $\alpha _{k}^{\ast }$,
$\beta _{k}^{\ast},k=1,...,k_{0}$ are defined by equality $(24)$.
Here the function $R_{012}(z)$ is holomorphic outside the set $\Omega $
and the points $\pm 1$.} \footnote{%
Let us remark that formula (27) does not include the case when one of
the points $\alpha _{k}^{*}$ or $\beta _{k}^{*}$ in equal to the infinity,
i.e. when $0$ belongs to the boundary of the set $V(\Omega _{1})$. In this
case the respective factor in $R_{0}(z)$ should be defined as
$P(z,\chi _{(\alpha_{k}^{*},\beta _{k}^{*})}n(k)\pi )$.}

\medskip

P r o o f. \quad The holomorphy of $R_{012}(z)$ outside the set
$\Omega \cup \{1,-1\}$ easily follows from its definition: the
function $R_{0}(z)R_{0}^{\prime }(z)R_{0}^{\prime \prime }(z)$
may have singularities only on the boundary of the set
$\Omega _{1}\cup \Omega _{2}\cup \{-1,1\}$, the function $R_{1}(z)$
may have singularities only on the set $\Omega _{1}$, and $R_{2}(z)$
on $\Omega _{2}$. According to (21), (25), (12),
$$
N_{1}(z)\ =\ CP(z,\nu )
$$
$$
   =CP(z,\Bigl(\chi _{1}+V(\chi _{1})+\chi _{2}+\chi _{\Delta}^{(k_{0})}
   +V(\chi _{\Delta }^{(k_{0})})+\chi _{\Delta }^{\prime }
   +V(\chi_{\Delta }^{\prime })+\chi _{\Delta }^{\prime \prime }
   +V(\chi _{\Delta}^{\prime \prime })\Bigr)\nu )=
$$
$$
   =CP(z,\chi _{1}\nu )P(z,V(\chi _{1})\nu )P(z,{\frac{1}{2}}\chi _{2}\nu )
    P(z,{\frac{1}{2}}\chi _{2}\nu )\times
$$
$$
   \times P(z,\chi_{\Delta }^{(k_{0})}\nu )
   P(z,V(\chi _{\Delta }^{(k_{0})}\nu)
   P(z,\chi _{\Delta }^{\prime }\nu )
   P(z,V(\chi _{\Delta }^{\prime })\nu)
   P(z,\chi _{\Delta }^{\prime \prime }\nu )
   P(z,V(\chi _{\Delta }^{\prime\prime })\nu ).\eqno(30)
$$
Here, according to the properties (14), (12), and the property $\nu (t)=-\nu
(t^{-1})$ (see (23)),
$$
   P(z,\chi_{1}\nu )P(z,V(\chi _{1})\nu)
   =P(z,\chi_{1}\nu )P(z^{-1},\chi_{1}\nu)=R_{1}(z)R_{1}(z^{-1}),
   \eqno(31)
$$
$$
   P(z,{\frac{1}{2}}\chi _{2}\nu )P(z,{\frac{1}{2}}\chi _{2}\nu )
   =P(z,{\frac{1}{2}}\chi _{2}\nu )P(z^{-1},{\frac{1}{2}}\chi _{2}\nu )=
     R_{2}(z)R_{2}(z^{-1}),  \eqno(32)
$$
$$
   P(z,\chi_{\Delta }^{\prime }\nu )P(z,V(\chi_{\Delta }^{\prime })\nu)
   =R_{0}^{\prime }(z)R_{0}^{\prime }(z^{-1}),\eqno(33)
$$
$$
   P(z,\chi _{\Delta }^{\prime \prime }\nu )
   P(z,V(\chi _{\Delta }^{\prime\prime })\nu )
   =R_{0}^{\prime \prime }(z)R_{0}^{\prime \prime }(z^{-1}),\eqno(34)
$$
and, from definition (4),
$$
   P(z,\chi _{\Delta }^{(k_{0})}\nu )P(z,V(\chi _{\Delta }^{(k_{0})})\nu )=
   \exp \Bigl\{{\frac{1}{{\pi }}}\int\limits_{-\infty }^{\infty }
  (\chi^{(k_{0})}(t)+V(\chi ^{(k_{0})})(t))\nu (t)
     ({\frac{1}{{t-z}}}-{\frac{t}{{1+t^{2}}}})dt\Bigr\}=
$$
$$
  =\prod_{k=1}^{k_{0}}\Bigl({\frac{{z-\beta _{k}}}{{z-\alpha _{k}}}}
   \Bigr)^{n(k)}\Bigl({\frac{{z-\alpha _{k}^{-1}}}{{z-\beta _{k}^{-1}}}}
   \Bigr)^{-n(k)}\exp \Bigl\{{\frac{1}{{\pi }}}\int\limits_{-\infty }^{\infty }
   (\chi^{(k_{0})}(t)+V(\chi ^{(k_{0})})(t))\nu (t){\frac{-t}{{1+t^{2}}}}dt\Bigr\}.
$$
The exponential in the r.h.s is a constant. Further, it is evident that the
right-hand side of the latest equality will not change if we replace
$\alpha_{k}$, $\varphi _{k}$ to $\alpha _{k}^{\ast }$,
$\varphi _{k}^{\ast }$. That is why
$$
P(z,\chi _{\Delta }^{(k_{0})}\nu )P(z,V(\chi ^{(k_{0})})_{\Delta }\nu)
  =C\prod_{k=1}^{k_{0}}\Bigl({\frac{{z-\beta _{k}^{\ast }}}
  {{z-\alpha_{k}^{\ast }}}}\Bigr)^{n(k)}
  \Bigl({\frac{{z-{\beta _{k}^{\ast }}^{-1}}}
  {{z-{\alpha _{k}^{\ast }}^{-1}}}}\Bigr)^{n(k)}.
$$
But
$$
   {\frac{{z-{\beta _{k}^{\ast }}^{-1}}}{{z-{\alpha _{k}^{\ast }}^{-1}}}}
   ={\frac{\alpha _{k}^{\ast }}{\beta _{k}^{\ast }}}\cdot {\frac{{z^{-1}-
   {\beta_{k}^{\ast }}}}{{z^{-1}-{\alpha _{k}^{\ast }}}}},
$$
so
$$
    P(z,\chi _{\Delta }^{(k_{0})}\nu )P(z,V(\chi _{\Delta }^{(k_{0})}\nu)=
    CR_{0}(z)R_{0}(z^{-1})\eqno(35)
$$
(with another constant $C>0$). Thus, from the equalities (30)--(35) follows
(26). \hfill \rule{0.5em}{0.5em}

\bigskip

Having factorized $N_{1}(z)$, we now factorize $N_{0}(z)$. We define
$$
R_{3}(z)=\hat{P}(z,{\frac{\hat{\nu}_{0}}{2}}).\eqno(36)
$$
For $|z|\neq 1$, according to (15), (16), taking into account the oddness of
$\hat{\nu}_{0}(\theta )$ and representation (21), we have
$$
R_{3}(z)R_{3}(z^{-1})=\hat{P}(z,{\frac{\hat{\nu}_{0}}{2}})\hat{P}(z^{-1},{%
\frac{\hat{\nu}_{0}}{2}})=\hat{P}(z,{\frac{\hat{\nu}_{0}}{2}})\hat{P}(z,{%
\frac{\hat{\nu}_{0}}{2}})=?N_{0}(z).\eqno(37)
$$

Let us define $R_\mu(z)$ by formula (7). It satisfies property (8). So, from
the definition of $R_3(z)$, the property (6) and the holomophy of
$R_{012}(z)$ on the unit circle (possibly, excepting the points $\pm 1$),
we have for the function
$$
R_{0123\mu}(z) \equiv R_0(z) R^{\prime}_0(z) R^{\prime\prime}_0(z)  R_1(z)
R_2(z) R_3(z) R_\mu(z) \eqno(38)
$$
the relation
$$
{\frac{ |R_{0123\mu}^+(\xi)| }{|R_{0123\mu}^-(\xi)| }} = {\frac{
|R_{012}^+(\xi)| }{|R_{012}^-(\xi)| }} \cdot {\frac{ |R_3^+(\xi)| }{%
|R_3^-(\xi)| }} \cdot {\frac{ |R_{\mu}^+(\xi)| }{|R_{\mu}^-(\xi)| }} =
\mu(\xi),  \quad |\xi|=1.
$$

But, as it follows directly from the definition of $R_\mu(z)$, with the use
of the evenness of $\mu(\theta)$,
$$
R_\mu(z) R_\mu(z^{-1}) = {\rm const} >0. \eqno(39)
$$
Thus, for the function $R_{0123\mu}(z)$, defined by (38), we have from (20),
(26), (39), (39)
$$
R_{0123\mu}(z)R_{0123\mu}(z^{-1}) = C N(z), \quad C>0.
$$
Finally, taking
$$
R(z) \equiv \sqrt{C} R_{0123\mu}(z),
$$
we see that $R(z)$ is the solution of our problem.

\hfill\rule{0.5em}{0.5em}\medskip

Certainly, the factorization (10) is not unique. As it will be seen from the
next theorem, with some additional restrictions to the function $N(z)$ we
can require additional conditions to the behavior of the function $R(z)$
near its singularities, for example, we can ask for the existence of the
limits in the metric of $L^{p},\ p\geq 1,$ of the functions $R(t\pm
i\varepsilon ),\ \varepsilon \rightarrow +0$.

{\bf Definition.} \ {\it Let $A$ be a certain set on the real axis, $%
U_{\delta }(A)$ be its $\delta $-neighborhood, and $f(z)$ be a holomorphic
function in $U_{\delta }(A)\backslash A$. We say that the function $f(z)$
locally belongs to the Hardy class $H^{p}$ in the neighborhood of the set
$A$, if for some $\delta >0$ the functions
$f(t\pm i\varepsilon ),\ t\in {\bf R}\cap U_{\delta }(A),$
converge as $\varepsilon \rightarrow +0$ \ in the metric of $L^{p}$. }

\medskip

{\bf Theorem 3.} \ {\it Let the function $N(z)$, the sets $\Sigma ,\Omega
_{1},\Omega _{2},V(\Omega _{1})$ and the function $\mu (\theta )$ satisfy
the conditions of the previous theorem. Assume that the set $\Omega _{2}$
can be covered with a finite number of mutually disjoint intervals $\delta
_{l}$, and on each of them
$$
{\rm ess}\sup_{t\in \delta _{l}}\arg N^{+}(t)-{\rm ess}\inf_{t\in \delta
_{l}}\arg N^{+}(t)<\pi \ .
$$
}

{\it Then the function $N(z)$ can be factored out so that $(10)$ and $(11)$
hold, and the function $R(z)$ locally depends to the Hardy class $H^{2}$ in
the neighborhood of the set $\Omega _{2}$. }

P r o o f. \ We will only show what changes should be done in the proof of
Theorem 1 to apply it to Theorem 2. The set
$$
\Delta ={\bf R}\backslash ([-1,1]\cup \Omega _{1}\cup V(\Omega _{1}))=\cup
_{k}\Delta _{k}\,,
$$
is introduced, where $\Delta _{k}=(\alpha _{k},\beta _{k})$ are mutually
disjoint intervals. The set $\Omega _{2}$ lie inside of these intervals. It
follows from the additional condition of the Theorem 3 that on every set $%
\Delta _{k}\backslash \Omega _{2}$ the function $n(t)=n(k)={\frac{1}{\pi }}%
\arg N(t)\in {\bf Z}$ is constant. We extend this function to the whole
interval $\Delta _{k}$ (possibly, including $\Omega _{2}$), and also we
extend it to $V(\Delta _{k})$ by the equality
$$
\tilde{n}(t)=%
\cases{n(k), \quad t \in \Delta_k, \cr
                     -n(k), \quad t \in V(\Delta_k). \cr}
$$
Further, we divide $\Delta _{k}$ to three groups: the intervals $\Delta
_{1},\ldots ,\Delta _{k_{0}}$, which have one of the endpoints in $V(\Omega
_{1})$ and the other in $\Omega _{1}\cup \{-1,1\}$, the intervals $\Delta
_{k}^{\prime }$, whose both endpoints belong to $V(\Omega _{1})$, and the
intervals $\Delta _{k}^{\prime \prime }$ (all the others). We choose numbers
$\alpha _{k}^{\ast }$, $\beta _{k}^{\ast }$, \ $k=1,\ldots ,k_{0}$, by the
same rule as it was done in the Theorem 2 (by formula (24). The functions $%
R_{0}(z)$, $R_{1}(z)$, $R_{3}(z)$, $R_{\mu }(z)$ are defined by formulae
(27), (29), (36), (7). We also define the functions
$$
R_{0}^{\prime }(z)\ =\ P(z,V(\chi _{\Delta }^{\prime })\tilde{n}\pi ),\qquad
R_{0}^{\prime \prime }(z)=P(z,\chi _{\Delta }^{\prime \prime }\tilde{n}\pi
),\qquad R_{2}(z)=P(z,{\frac{1}{2}}\chi _{2}(\nu -\tilde{n}\pi )).
$$
The rest of the proof does not change. We will explain how to prove that the
function $R(z)$ locally belongs to the Hardy space in the neighborhood of $%
\Omega _{2}$. The function $R_{2}(z)$ is represented in the multiplicative
form $R_{2}(z)=P(z,{\frac{1}{2}}\chi _{2}(\nu -\tilde{n}\pi ))$. Here, the
oscillation of the function ${\frac{1}{2}}\chi _{2}(t)(\nu (t)-\tilde{n}%
(t)\pi )$ on the interval $\delta _{l}$ is less than ${\frac{\pi }{2}}$.
This implies (see, e.g., [5])\footnote{%
We use the following fact from this book \newline
{\sl Theorem} \ (Smirnov). \ {\it If $f(z)$ is analytic in the disk $|z|<1$
and ${\rm Re}\,f(z)\geq 0$, than $f\in H^{p}$ for every $p>1$.} \newline
(Here the Hardy space in the unit disk $|z|<1$ is denoted by $H^{p}$).},
that $R_{2}(z)$ locally belongs to the Hardy space in the neighborhood of
any compact subset of the interval $\delta _{l}$. (Such intervals $\delta
_{l}$ cover $\Omega _{2}$ according to the condition of the theorem).
Simultaneously, the functions $R_{0}(z)$, $R_{0}^{\prime }(z)$, $%
R_{0}^{\prime \prime }(z)$, $R_{1}(z)$, $R_{3}(z)$ and $R_{\mu }(z)$ are
holomorphic on $\Omega _{2}$, from where we obtain that $R(z)$ locally
belongs to the Hardy class $H^{2}$ in the neighborhood of the set $\Omega
_{2}$.

\hfill\rule{0.5em}{0.5em}\medskip

We remark that the functions $R_{0}^{\prime }(z)$, $R_{0}^{\prime \prime }(z)
$ can also be presented in the form of a product (27) (generally speaking,
this product is infinite). However, if the intervals $\Delta _{k}^{\prime
\prime }$ do not belong to a finite set of the real axis, then additional
factors will occur in the expression for $R_{0}^{\prime \prime }(z)$.

\bigskip

{\bf Remark.} This paper is the translation from [6].

\bigskip
\bigskip
\bigskip

\centerline{\bf References}

\bigskip

\begin{flushleft}

[1] {\it S.P.~Novikov ed.,} \ Soliton theory. Inverse problem method.
\ Moscow, 1980\\

\smallskip

[2] {\it Anne Boutet de Monvel and Vladimir Marchenko,} \
The Cauchy problem for nonlinear Schr\"odinger equation with bounded
initial data. \ Mat. fiz., analiz, geom. (1997),
v.4, ü~1/2, p.3--45.\\

\smallskip

[3] {\it P. Koosis,} \ Introduction in the theory of spaces $H^p$. \
Moscow (1984).\\

\smallskip

[4] {\it N.I.~Akhiezer and I.M.~Glazman,} \ Theory of linear
operators in Hulbert space. \ Pitman, Boston (1981)\\

\smallskip

[5] {\it J.V. Garnett,} \ Bounded analytic functions. \
Academic press, New York, London (1981).

\smallskip

[6] {\it M.Kudryavtsev}, \
The Riemann problem with additional singularities {\sl (in Russian)}.
Mat. fiz., analiz, geom. (2000), v.~7, ü~2, p.~196--208.

\end{flushleft}

\end{document}